\documentclass[10pt]{article}
\usepackage[top=1in,bottom=1.1in,left=1.25in,right=1.25in]{geometry}

\usepackage{multirow}
\usepackage{color, colortbl}
\usepackage{caption}
\usepackage{subcaption}
\usepackage{orcidlink}

\usepackage{graphicx}
\usepackage{ragged2e}
\usepackage{tabto}
\usepackage{amssymb}
\usepackage{float}
\usepackage{amsmath}
\usepackage{amsthm}
\usepackage{epsfig}
\usepackage{amscd}
\usepackage{color}
\usepackage{url}

\usepackage{hyperref}
\hypersetup{
    colorlinks=true,
    linkcolor=black,
    citecolor=blue,
    urlcolor=cyan,
    pdfpagemode=FullScreen,
    }

\usepackage[sort,compress,square,numbers]{natbib}

\usepackage{cleveref}
\usepackage{afterpage}
\usepackage{fancybox}
\usepackage{ifthen}
\usepackage{multicol}
\usepackage[font=footnotesize]{caption}
\usepackage{bigints}

\usepackage[sf,bf,medium]{titlesec}
\usepackage{setspace}

\usepackage{dsfont}
\usepackage{isomath}
\usepackage{doi}
\newcommand{\mtx}[1]{\mathrm{#1}}

\floatstyle{plain}
\restylefloat{figure}

\providecommand{\keywords}[1]
{
  \small	
  \textbf{\textit{Keywords---}} #1
}

\newcommand\cT{\mathcal T}

\newcommand\cO{\mathcal O}

\newcommand\bc{\boldsymbol c}

\newcommand\bn{\boldsymbol n}

\newcommand\bx{\boldsymbol x}

\newcommand\by{\boldsymbol y}

\newcommand\dmin{d_{\textrm{min}}}

\newcommand\gammatop{\Gamma \rightarrow P}
\newcommand\ptogamma{P \rightarrow \Gamma}

\newcommand\cA{\mathcal{A}}

\newcommand\cS{\mathcal{S}}

\newcommand\cD{\mathcal{D}}

\newcommand\bbC{\mathbb C}

\newcommand\bbR{\mathbb R}

\newcommand\efmm{\varepsilon_{\textrm{mv}}}
\newcommand\egmres{\varepsilon_{\textrm{gmres}}}
\newcommand\timeit{t_{\textrm{it}}}
\newcommand\timetot{t_{\textrm{tot}}}
\newcommand\nit{n_{\textrm{it}}}

\newcommand\uin{u^{\textrm{in}}}
\newcommand\dudnin{\frac{\partial u}{\partial n}^{\textrm{in}}}
\newcommand\usc{u^{\textrm{sc}}}
\newcommand\dudnsc{\frac{\partial u}{\partial n}^{\textrm{sc}}}
\theoremstyle{definition}

\theoremstyle{remark}
\newtheorem{remark}{Remark}

\numberwithin{equation}{section}

\crefname{equation}{}{}

\begin{document}

%\maketitle

\begin{titlepage}

  \raggedleft
  % {\small \texttt{CCM Technical Report\\
  {\small \texttt{STATUS: arXiv pre-print\\
    \today}}
  
  \hrulefill

  \vspace{4\baselineskip}

  \raggedright {\Large \sffamily\bfseries On the construction of
    scattering matrices for irregular\\ \vspace{.25\baselineskip}
    or elongated enclosures using Green's representation formula}
  
  \vspace{5\baselineskip}
  % author 1
  \normalsize Carlos Borges\footnote{Research supported in part by the Office of Naval Research under award N00014-21-1-2389.}\,\orcidlink{0000-0002-8236-6131} \\
  \small \emph{University of Central Florida, Department of Mathematics\\
    Orlanda, FL 32816}
 
  \texttt{carlos.borges@ucf.edu}
  
 \vspace{1.25\baselineskip}
  \normalsize Leslie Greengard\,\orcidlink{0000-0003-2895-8715} \\
  \small \emph{Center for Computational Mathematics\\
    Flatiron Institute  \\
    New York, NY 10010} and \\
   \emph{Courant Institute, New York University \\
    New York, NY 10012} \\
  \texttt{lgreengard@flatironinstitute.org}

 \vspace{1.25\baselineskip}
    \normalsize Michael O'Neil\,\orcidlink{0000-0003-2724-215X} \\
    \small \emph{Center for Computational
    Mathematics\\ Flatiron Institute\\
    New York, NY 10010}\\
    \texttt{moneil@flatironinstitute.org}
  
 \vspace{1.25\baselineskip}
    \normalsize Manas Rachh\,\orcidlink{0000-0001-5488-3839}\\
    \small \emph{Center for Computational
    Mathematics\\ Flatiron Institute\\
    New York, NY 10010}\\
    \texttt{mrachh@flatironinstitute.org}

\end{titlepage}

\begin{abstract}%must be between150-250 words for JSC
  Multiple scattering methods are widely used to reduce the computational
  complexity of acoustic or electromagnetic scattering problems when waves
  propagate through media containing many identical inclusions. Historically,
  this numerical technique has been limited to situations in which the
  inclusions (particles) can be covered by nonoverlapping disks in two
  dimensions or spheres in three dimensions. This allows for the use of
  separation of variables in cylindrical or spherical coordinates to represent
  the solution to the governing partial differential equation. Here, we provide
  a more flexible approach, applicable to a much larger class of geometries. We
  use a Green's representation formula and the associated layer potentials to
  construct incoming and outgoing solutions on rectangular enclosures.  The
  performance and flexibility of the resulting scattering operator formulation
  in two-dimensions is demonstrated via several numerical examples for
  multi-particle scattering in free space as well as in layered media. The
  mathematical formalism extends directly to the three dimensional case as well,
  and can easily be coupled with several commercial numerical PDE software
  packages.
\end{abstract}

\keywords{Helmholtz equation, multiple particle scattering, scattering matrix, irregular or elongated enclosures} 

\tableofcontents
\newpage

\normalsize
\setstretch{1.1}

\section{Introduction}

A powerful tool in the analysis of wave propagation problems in domains with
many inclusions is multiple scattering theory (see, for example,
\cite{botten2001photonic,FMPS,GroteKirsch,Hesford,Ioannidou,FMPS2D,laili19,Xu}
and the monographs \cite{Bohren,Martin,Hafner,MTL}).  The workhorse in this
approach is the construction of the {\em scattering operator} or \emph{matrix}
for an irregular compact inclusion -- typically obtained from computing the map
from incoming data to outgoing scattered solutions on an enclosing disk (circle)
in the plane or an enclosing ball (sphere) in three dimensions.  We will refer
to the enclosing surfaces as {\em proxy surfaces}. For simplicity, here we
restrict our attention to the two-dimensional acoustic scattering case
(extensions to three dimensions and other elliptic PDEs are relatively
straightforward), where the governing equation is the scalar Helmholtz equation
\begin{equation} 
\Delta u + k^2 u = 0, 
\label{helmeq}
\end{equation}
with~$k$ a complex wavenumber with non-negative imaginary part.
Given a disk $D$ of radius $R$, 
solutions of \eqref{helmeq} that are regular inside the disk can be represented as
\begin{equation} 
\sum_{n = - \infty}^\infty \alpha_n \, J_n(k r) \, e^{i n \theta},
\label{jexp}
\end{equation}
where~$(r,\theta)$ are the polar coordinates of a point with respect to the disk
center and $J_n$ denotes the $n$th-order Bessel function of the first kind.
Solutions of \eqref{helmeq} that are regular outside the disk and which satisfy
the Sommerfeld radiation condition
\begin{equation}
\lim_{r \rightarrow \infty}
r^{1/2} \left( \frac{\partial u}{\partial r} - i k u \right) = 0
\label{sommerfeld}
\end{equation}
can be represented as
\begin{equation} 
\sum_{n = - \infty}^\infty \beta_n \, H_n(k r) \, e^{i n \theta},
\label{hexp}
\end{equation}
where $H_n$ is the $n$th-order Hankel function of the first kind.  In this basis
of Bessel functions, given an incoming field sampled on the boundary of the disk
as a finite Fourier series, the scattering operator maps the vector of Fourier
coefficients~$\{ \alpha_n \}$ to the corresponding coefficients $\{ \beta_n \}$
of the truncated outgoing expansion~\eqref{hexp}.  The {\em construction} of the
scattering operator, of course, requires the one-time solution of a collection
of boundary value problems on the inclusion itself, the details of which depend
on the material properties of the inclusion.

\begin{figure}[!t]
  \centering
  \includegraphics[width=.75\linewidth]{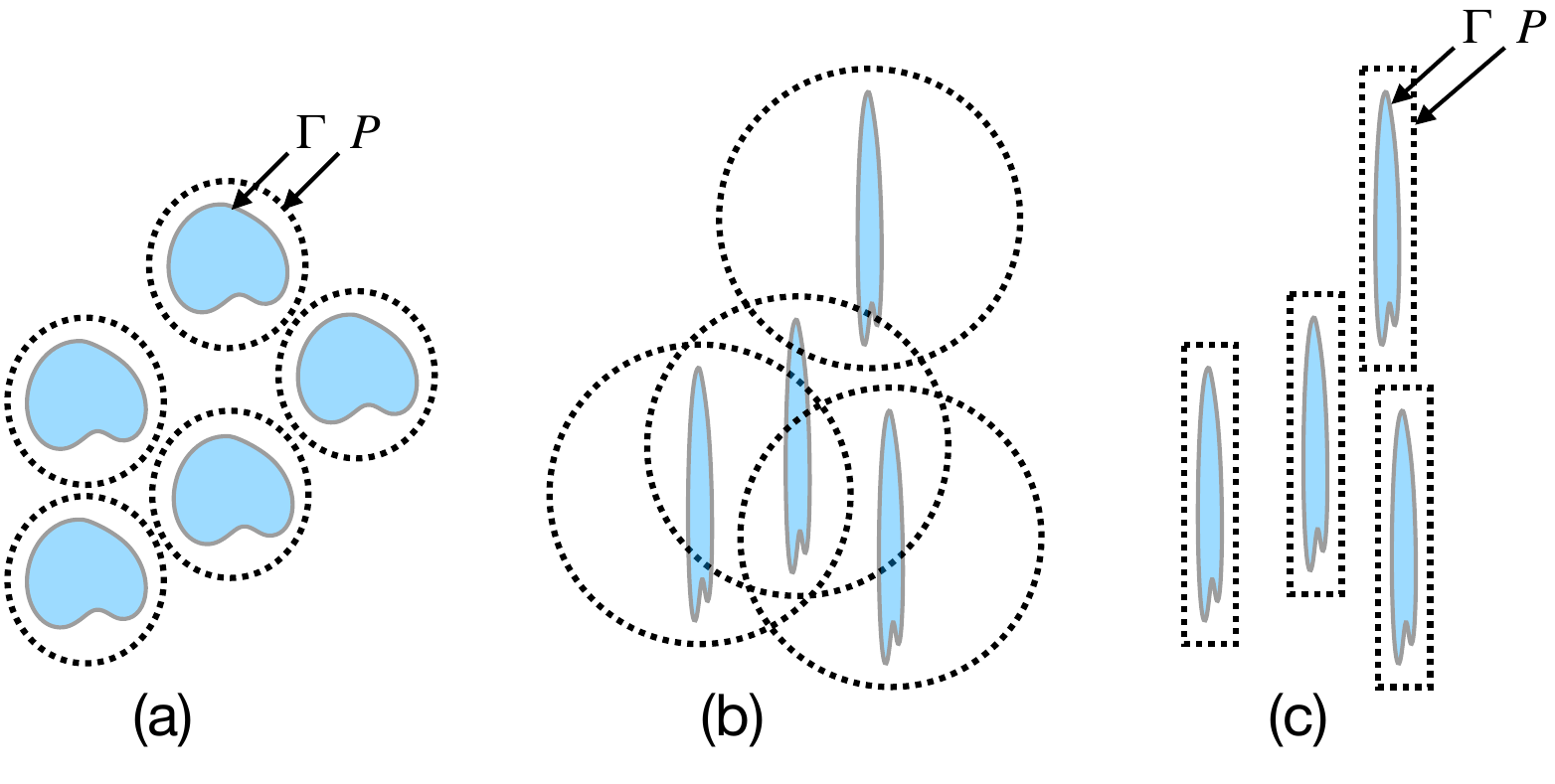}
  \caption{In multiple scattering methods, each inclusion with
    boundary~$\Gamma$ is enclosed in a proxy surface~$P$. For any incoming
    field, sampled on~$P$, the scattering matrix yields a representation for the
    outgoing field induced on~$\Gamma$ (but only in the {\em exterior}
    of~$P$). These methods work well when inclusions are spatially separated.
    When~$P$ is a disk, as in (a), even fairly close interactions can be
    accurately computed, requiring only that the enclosing disks are
    disjoint. For high aspect ratio inclusions, as in (b), they fail since the
    enclosing disks overlap even when there is a significant gap between them.
    Here, we extend the approach by constructing scattering matrices on more
    closely contoured enclosures such as a rectangle, as depicted in (c).}
  \label{scatmatfig}
\end{figure}

This classical theory is simple because the expansions~\eqref{jexp}
and~\eqref{hexp} are explicit bases for solutions of the Helmholtz equation
(obtained by separation of variables).  In practice, however, this approach
fails when the inclusions are moderately close to touching and far from circular
in shape (see Fig.~\ref{scatmatfig}) due to the slow convergence 
or divergence of such expansions.

In this paper, we describe a new technique for acoustic or electromagnetic
scattering from complex microstructured materials in two dimensions, governed
by~\eqref{helmeq}, assuming only that the inclusions are identical up to
rotation and translation (or drawn from a small set of such particles).  The
power of the approach, as with classical multiple scattering methods, is that it
yields efficient algorithms for microstructure design -- in the context of an
outer optimization process, each new configuration requires the calculation of
an objective function which, in general, depends on the scattered field. With
thousands of inclusions, this is intractable when done directly. The use of
scattering matrices can dramatically reduce the number of degrees of freedom in
such calculations as well as the condition number of the full problem so that
such design loops become computationally feasible, especially when coupled with
fast algorithms to compute the interactions between the individual inclusions.

We summarize the necessary background from potential theory in the next section,
followed by a description of the construction of generalized scattering
operators. We then apply the method to rectangular enclosures and demonstrate its
performance in both free space and layered media.  We will use the terms
\emph{enclosure} and \emph{proxy surface} interchangeably.  It is worth noting
that the solver used for modeling individual inclusions when constructing the
scattering matrix can be called in ``black box'' fashion, and easily coupled to
any commercial or non-commercial software packages.

\begin{remark}
  There is an alternative to defining scattering matrices on enclosures that has
  emerged over the last decade or so.  Namely, modern fast direct solvers can
  (a) sample the incoming field on the obstacle itself, and (b) solve an
  associated integral equation to obtain an equivalent charge or dipole
  distribution along the surface of the obstacle (or an enclosing proxy surface)
  that accurately reproduces the
  scattered solution in the far field.  This representation can be compressed by
  identifying a subset of points on the surface that are sufficient to represent
  both the incoming and scattered fields.  The map from the incoming field {\em
    sampled} at those points to the charge/dipole strengths {\em induced} at
  those points is itself a kind of scattering matrix.  This idea was used in a
  principled hierarchical manner to construct fast direct solvers for
  non-oscillatory elliptic integral equations on planar domains with corners
  in~\cite{bremer2012corner} and on smooth surfaces in three dimensions
  in~\cite{bremer-quad-2015}. We will refer to this approach as
  {\em skeletonization}. In some applications it may be beneficial to
  construct scattering matrices through skeletonization, while in others
  it may be more practical to build a few
  ``standard'' scattering matrices in the sense depicted in Fig.~\ref{scatmatfig}. We will return to this question 
  in the concluding section.
\end{remark}

\section{Mathematical preliminaries}

In this section, we summarize the properties of layer potentials using the
Green's function for the scalar Helmholtz equation \eqref{helmeq}, given by
\begin{equation}
g_{k}(\bx,\by) = \frac{i}{4} H_{0}^{(1)}(k |\bx - \by|) \, ,
\end{equation}
where $H_0^{(1)}$ is the zeroth-order Hankel function of the first kind.  The
variables~\mbox{$\bx = (x_1,x_2)$} and~\mbox{$\by=(y_1,y_2)$} denote points
in~$\bbR^2$.  See~\cite{ColtonKress} for a thorough discussion of this Green's
function.

\subsection{Layer potentials}

Given a closed curve~$\Gamma$ and densities~$\sigma$ and~$\mu$ supported along
this curve, we define the standard single and double layer
potentials~$\cS_{\Gamma} [\sigma]$ and~$\cD_{\Gamma}[\mu]$ by
\begin{equation}
  \begin{aligned}
    \cS_{\Gamma}[\sigma](\bx) &= \int_{\Gamma} g_{k}(\bx,\by) \, \sigma(\by) \, ds, \\
    \cD_{\Gamma}[\mu](\bx) &= \int_{\Gamma} \left( \bn(\by) \cdot \nabla_{\by}
                                g_{k}(\bx,\by)  \right)  \mu(\by) \, ds ,
  \end{aligned}
\end{equation}
where it is assumed that~$\bx \notin \Gamma$.  These potentials automatically
satisfy the Sommerfeld radiation condition~\eqref{sommerfeld}, and furthermore,
for $\bx \notin \Gamma$ they are infinitely differentiable functions. For
values of~$\bx$ on the boundary~$\Gamma$, the potentials are well-defined in
terms of weakly singular integral operators, denoted by~$S_\Gamma$
and~$D_\Gamma$.  For~$\bx_0 \in \Gamma$, these satisfy the jump relations
\begin{equation}
  \label{eq:jumprel}
  \begin{aligned}
    \lim_{\substack{\bx \to  \bx_{0} \\ \bx \in \Omega^{\pm}}} \cS_{\Gamma} [\sigma](\bx)
    &= \mtx{S}_{\Gamma}[\sigma](\bx_{0}) \, ,\\
    \lim_{\substack{\bx \to  \bx_{0} \\ \bx \in \Omega^{\pm}}} \cD_{\Gamma} [\mu](\bx)
    &= \pm \frac{\mu(\bx_{0})}{2} + \mtx{D}_{\Gamma}[\mu](\bx_{0}) \, .
%\lim_{\substack{\bx \to  \bx_{0} \\ \bx \in \Omega^{\pm}}} \cS'_{\Gamma} [\sigma] &= \pm \frac{\sigma(\bx_{0})}{2} + \mtx{S}'_{\Gamma}[\sigma](\bx_{0}) \, ,\\
%\lim_{\substack{\bx \to  \bx_{0} \\ \bx \in \Omega^{\pm}}} \cD'_{\Gamma} [\sigma] &= \mtx{D}'_{\Gamma}[\sigma](\bx_{0}) \, .
\end{aligned}
\end{equation}
Next, consider an enclosure or proxy surface~$P$ with $\Gamma$ in its interior,
as shown earlier in Fig.~\ref{scatmatfig}. Let the unit outward normal
along~$P$ be denoted by~$\bn$.  The normal derivatives
of~$\cS_{\Gamma} [\sigma]$ and~$\cD_{\Gamma}[\mu]$ along~$P$ are defined as:
\begin{equation}
  \begin{aligned}
\cS'_{\gammatop}[\sigma](\bx) &= \bn(\bx) \cdot \nabla
\cS_{\Gamma} [\sigma](\bx) ,\\
\cD'_{\gammatop}[\mu](\bx) &= \bn(\bx) \cdot \nabla_{\bx} 
\cD_{\Gamma} [\mu](\bx)   .
\end{aligned}
\end{equation}
In general, given two disjoint oriented curves~$\Gamma_{1}$ and~$\Gamma_{2}$, we
will denote by~$\cS_{\Gamma_{1} \rightarrow \Gamma_{2}}[\sigma]$ the restriction
of the layer potential~$\cS_{\Gamma_{1}}[\sigma]$ to points~$\bx \in
\Gamma_2$. (The same notation will be used for $\cD$.)  In a slight abuse of
notation, when $\Gamma_2=\Gamma_1=\Gamma$, we will refer to the weakly singular
operators~$S_{\Gamma \rightarrow \Gamma}$ and~$D_{\Gamma \rightarrow \Gamma}$ as
$S_{\Gamma}$, $D_\Gamma$ (as was done above).  Finally, we will denote a function
$u: \bbR^2 \to \bbC$ restricted to a closed curve~$P$ by~$u|_P$.

\subsection{Combined field integral equation}

While the scattering operator construction described below is independent of the
material properties of the inclusion, for the sake of concreteness, we assume
that, using the language of acoustics, it is
\emph{sound-soft}~\cite{ColtonKress}. That is, given an incoming (pressure)
field $\uin$, the total pressure $u = \uin + \usc$ must satisfy the homogeneous
Dirichlet condition $u=0$ so that $\usc = -\uin$ on the inclusion boundary
$\Gamma$.  The classical integral equation for such problems makes use of the
\emph{combined field} representation
\begin{equation}
  \usc(\bx) = \cD_{\Gamma} [\sigma](\bx) + ik \cS_{\Gamma}[\sigma](\bx) \,,
  \quad \bx \in \Omega^{c} \, .
\end{equation}
Using the jump relations in~\cref{eq:jumprel}, $\sigma$ can be shown to be given
by \cite{ColtonKress}
\begin{equation}
\label{eq:sig_sol}
\sigma = -\left(\frac{1}{2} + 
\mtx{D}_{\Gamma}+ik \mtx{S}_{\Gamma}\right)^{-1} \uin|_{\Gamma} \, .
\end{equation}

\subsection{Green's identities}

Suppose that a function $u$ satisfies the Helmholtz equation in the
interior~$\Omega$ of a closed curve~$\Gamma$. Green's representation formula
\cite{ColtonKress} states that for any~$\bx \in \Omega$,
\begin{equation}
  u(\bx) = -\cD_\Gamma[u](\bx) + 
  \cS_\Gamma \left[ \frac{\partial u}{\partial n} \right](\bx)
  \label{greenrepint}
\end{equation}
where~$\partial/\partial n$ denotes the directional derivative in the direction
of~$\bn$, the unit \emph{outward} normal to~$\Gamma$ (i.e. the direction pointing
into the unbounded domain~$\Omega^c$).  Furthermore, assuming a function~$v$
satisfies the Helmholtz equation in the region~$\Omega^{c}$ exterior
to~$\Omega$, we have
\begin{equation}
  v(\bx) = \cD_\Gamma[v](\bx) - 
  \cS_\Gamma \left[ \frac{\partial v}{\partial n} \right](\bx).
  \label{greenrepext}
\end{equation}
The flip in signs is due to the definition of the normal direction~$\bn$.

\section{Scattering operator construction}
\label{sec:scat}

We now turn to describing a method for constructing a scattering operator using
nothing more than the Green's representation identities given in the previous
section.

Let $P$ be a rectangle that encloses a sound-soft obstacle~$\Omega$ with
boundary $\Gamma$, and let $\Omega_{P}$ denote the interior of $P$.  The
scattering operator in this formulation is a map from values of the incoming data
$(\uin, \dudnin)$ on $P$ to samples of the scattered field $(\usc,\dudnsc)$ on
$P$. Since $\uin$ is assumed to satisfy the homogeneous Helmholtz equation
inside $P$, we may use Green's representation formula \eqref{greenrepint} to
obtain
\begin{equation}
\label{eq:uin_green}
\uin(\bx) = -\cD_{P}\left[\uin\right](\bx)
+ \cS_{P}\left[\dudnin\right](\bx) \,, \quad \bx \in \Omega_{P} \, .
\end{equation}
This is a parameterization of~$\uin$ in terms of its boundary data
along~$P$. Given~$\uin$ along~$\Gamma$, we can solve~\eqref{eq:sig_sol} to
obtain the density~$\sigma$ along~$\Gamma$ and then evaluate
the scattered field for $\bx \in P$ as
\begin{equation}
  \label{eq:usc-rep}
  \begin{aligned}
    \usc(\bx) &= \cD_{\gammatop}[\sigma](\bx) + ik \cS_{\gammatop}[\sigma](\bx) \,, \\ 
    \dudnsc(\bx) &= {\cD}'_{\gammatop}[\sigma](\bx)
                   + ik {\cS}'_{\gammatop}[\sigma](\bx) \, .
  \end{aligned}
\end{equation}
Finally, combining equations~\cref{eq:uin_green,eq:sig_sol,eq:usc-rep}, the
scattering operator which maps the incoming data and its normal derivative to
the scattered field and its normal derivative, which we will denote
by~$\mtx{A}^{\Gamma}_{P}$, is given by the following composition of operators:
\begin{equation}
\label{eq:scatmat}
  \mtx{A}^{\Gamma}_{P}
  = 
  \begin{bmatrix}
    {\cD}_{\gammatop} + ik {\cS}_{\gammatop} \\
    {\cD}'_{\gammatop} + ik {\cS}'_{\gammatop}
  \end{bmatrix}
  \left( \frac{1}{2} + \mtx{D}_{\Gamma} + ik\mtx{S}_{\Gamma} \right)^{-1}
  \begin{bmatrix}
    {\cD}_{\ptogamma} & -{\cS}_{\ptogamma}
  \end{bmatrix} \, .
\end{equation}
That is to say,
\begin{equation}
\begin{bmatrix}
    \usc|_{P} \\
    \dudnsc|_{P}
\end{bmatrix}
= \mtx{A}^{\Gamma}_{P} 
\begin{bmatrix}
\uin|_{P}  \\
\dudnin|_{P} 
\end{bmatrix} \, .
\end{equation}

Lastly, suppose that $\Gamma$ and $\Gamma'$ are related by a unitary affine
transformation, and that the corresponding proxy surfaces $P$ and $P'$ are
related by the same unitary affine transformation. Then, due to the translation
and rotation invariance of the Helmholtz operator, the associated scattering
operators are identical,
i.e.~$\mtx{A}^{\Gamma}_{P} = \mtx{A}^{\Gamma'}_{P'}$.

\subsection{Multi-particle scattering}
\label{sec:multi}

We now consider the classical problem of scattering from a collection of~$M$
obstacles~$\Omega_j$ with boundaries given by~$\Gamma_{j} = \partial\Omega_j$,
for~$j = 1,\dots,M$. As above, let~$P_{j}$ denote a rectangle
enclosing~$\Gamma_{j}$ and let $\Omega_{P_{j}}$ denote its interior. In
general, and in what follows, there is no assumption that the obstacles are
identical or related by affine transformations. If this happens to be the case,
then each~$P_j$ can be constructed from an affine transformation of a
canonical~$P$, and the resulting scattering operators will be identical.

We may represent the total field in the exterior of all proxy surfaces by
\begin{equation}
\label{mpsfree}
u = \sum_{j=1}^{M} \left(\cD_{P_{j}}\left[\usc \right] - \cS_{P_{j}} \left[\dudnsc \right] \right) + \uin \, ,
\end{equation}
using Green's representation theorem.  Since the scattering operator accounts
for the response of each inclusion to an arbitrary incoming field, the condition
that needs to be satisfied is simply the continuity of the potential and its
normal derivative on the (artificial) proxy surfaces boundaries.  That is, on
enclosure $P_i$, we must have
\begin{multline}
\begin{bmatrix}
    \usc|_{P_{i}} \\
    \dudnsc|_{P_{i}}
\end{bmatrix}
- \sum_{\substack{j=1\\j \neq i}}^{M} 
\begin{bmatrix}
{\cD}_{P_{j} \rightarrow P_{i}} & -{\cS}_{P_{j} \rightarrow P_{i}} \\
{\cD}'_{P_{j} \rightarrow  P_{i}} & -{\cS}'_{P_{j} \rightarrow P_{i}} 
\end{bmatrix}\begin{bmatrix}
    \usc|_{P_{j}} \\
    \dudnsc|_{P_{j}}
\end{bmatrix} \\
 = 
 \mtx{A}^{\Gamma}_{P_{i}} 
\left(\begin{bmatrix}
\uin|_{P_{i}}  \\
\dudnin|_{P_{i}} 
\end{bmatrix} + \sum_{\substack{j=1\\j \neq i}}^{M} 
\begin{bmatrix}
{\cD}_{P_{j} \rightarrow P_{i}} & -{\cS}_{P_{j} \rightarrow P_{i}} \\
{\cD}'_{P_{j} \rightarrow  P_{i}} & -{\cS}'_{P_{j} \rightarrow P_{i}} 
\end{bmatrix}\begin{bmatrix}
    \usc|_{P_{j}} \\
    \dudnsc|_{P_{j}}
\end{bmatrix}\right) \, .
\label{mpsfreesys}
\end{multline}
Omitting the algebra, the left-hand side of \eqref{mpsfreesys} is the
contribution of the scattered field from scatterer~$\Omega_{i}$ to the net
scattered field, as seen from representation~\eqref{mpsfree}, and the right hand
side encodes the response of the $i$th inclusion to the total incoming field.
The above equations provide the setup for the classical multiple scattering
formalism~\cite{Martin,FMPS}.

Let the operators~$\cA$,~$\mtx{T}_{ij}$, and~$\cT$ be given by
\begin{equation}
\cA = 
\begin{bmatrix}
\mtx{A}_{P_{1}}^{\Gamma_{1}} & 0 & \ldots & 0 \\
0 & \mtx{A}_{P_{1}}^{\Gamma_{2}}  & \ldots & 0 \\
\vdots & \vdots & \ddots & \vdots \\
0 & 0 & \ldots & \mtx{A}_{P_{M}}^{\Gamma_{M}}
\end{bmatrix} \, ,
\label{Adef}
\end{equation}
\begin{equation}
T_{ij} = \begin{bmatrix}
{\cD}_{P_{j} \rightarrow P_{i}} & -{\cS}_{P_{j} \rightarrow P_{i}} \\
{\cD}'_{P_{j} \rightarrow P_{i}} & -{\cS}'_{P_{j} \rightarrow P_{i}} 
\end{bmatrix} \,.
\label{Tijdef}
\end{equation}
and 
\begin{equation}
\cT = 
\begin{bmatrix}
0 & \mtx{T}_{12} & \ldots &\mtx{T}_{1M} \\
\mtx{T}_{21} & 0 & \ldots & \mtx{T}_{2M} \\
\vdots & \vdots & \ddots & \vdots \\
\mtx{T}_{M1} & \mtx{T}_{M2} &\ldots & 0
\end{bmatrix} \, ,
\label{Tdef}
\end{equation}
Using these definitions,
the classical multiple scattering system~\eqref{mpsfreesys}
can be written concisely in linear algebraic form as:
\begin{equation}
\label{eq:multi-part-scat-main}
\left(I - (\cA+I) \cT \right) \begin{bmatrix}
\usc|_{P_{1}} \\
\dudnsc|_{P_{1}} \\
\vdots \\
\usc|_{P_{M}} \\
\dudnsc|_{P_{M}} 
\end{bmatrix}
= \cA
\begin{bmatrix}
\uin|_{P_{1}} \\
\dudnin|_{P_{1}} \\
\vdots \\
\uin|_{P_{M}} \\
\dudnin|_{P_{M}} 
\end{bmatrix} \, .
\end{equation}
The right hand side above only needs to be computed once, and assuming that the
inclusions are not extremely close-to-touching, solving the above system using
an iterative solver, such as GMRES, usually only requires a modest number of
iterations (particularly when a diagonal preconditioner is used~\cite{FMPS}).

\section{Multiple scattering in a layered medium}
\label{sec:scat-half}

In many applications, including meta-material
design~\cite{yu2014flat,blankrot2019efficient}, a collection of
identically-shaped scatterers are embedded in a layered medium to achieve some
response of interest (such as focusing an incoming wave).  For this, the
multiple scattering framework above needs to be modified to include the
effect of the layered medium itself.  This has been done previously using
classical scattering matrices (see, for example, \cite{FMPS2D,lim1992multiple}).

We show in this section how to extend our approach in the simplest setting: when
a collection of $M$ obstacles with boundaries~$\Gamma_{j}$ as well as their
enclosures $P_j$, $j =1,2,\ldots M$, lie in the upper half
plane~$\Omega_{+} = \{ (x_1,x_2)\,| \, x_{2} > 0 \}$.  We denote by
$\Omega_{-} = \{ (x_1,x_2)\,| \, x_{2}<0 \}$ the lower half plane and assume
that~$x_2=0$ is an {\em acoustic interface} with wavenumbers in the upper and
lower half planes given by $k_{\pm}$, respectively. A canonical transmission
condition along~$x_2=0$ is that the total field and its normal derivative are
continuous across the interface.  We let~$\usc_{\pm}$ denote the scattered field
in the upper or lower half-space, respectively, and with a slight abuse of
notation, will use~$\usc$ without the subscript when the context is clear.  The
extension to layered media is somewhat detailed but straightforward, so long as
each inclusion is completely contained in a single layer.

\subsection{The layered medium Green's function}

The simplest way to account for the presence of the infinite acoustic interface
is through the construction of the layered medium Green's function 
$g^{LM}$.  In what follows, we will assume that the \emph{source}~$\by$ is located in the upper half plane and that the target~$\bx$ can be located anywhere in~$\bbR^2$. An analogous definition of~$g^{LM}$ can be constructed when the source is in the lower half plane. The layered medium Green's function is then defined as
\begin{equation}
g^{LM}(\bx, \by) = 
\begin{cases}
g_{k_+}(\bx,\by) + s_{+}(\bx, \by), \quad &x_{2} > 0\\ 
s_{-}(\bx, \by), \quad &x_{2} < 0
\end{cases}
\end{equation}
where the corrections~$s_{\pm}$ are
outgoing solutions to the Helmholtz equation which enforce the 
continuity of~$g^{LM}$ along $x_2=0$:
\begin{equation}
  \begin{aligned}
(\Delta + k_{\pm}^{2}) s_{\pm} &= 0\,, &\quad \bx &\in \Omega_{\pm} \, , \\
s_{+} - s_{-} &= - g_{k_+}, 
&  \text{for } x_2 &= 0 ,\\
\frac{\partial s_{+}}{\partial x_{2}} - \frac{\partial s_{-}}{\partial x_{2}}  
&= - \frac{\partial g_{k_+} }{\partial x_2} \, ,
&  \text{for } x_2 &= 0.
\label{continuityx2}
\end{aligned}
\end{equation}
A discussion of the outgoing radiation conditions for solutions to the Helmholtz
equation in layered media can be found~\cite{chandler2007,epstein2023solving1,epstein2023solving2,epstein2024solving3}. It
is well-known that by applying the Fourier transform along the interface,
solutions (i.e. these corrections~$s_\pm$) to the homogeneous Helmholtz equation
can be represented as Sommerfeld integrals~\cite{chew,sommerfeld} of the form

\begin{equation}
s_{\pm}(\bx) 
= \int_{-\infty}^{\infty} 
  e^{-\sqrt{\xi^2 - k_{\pm}^2}|x_2|} \left(
\frac{\hat{\sigma}(\xi)}{\sqrt{\xi^2 - k_{\pm}^2}} \pm \hat{\mu}(\xi) \right) \, e^{i x_{1} \xi  } \, d\xi ,
\end{equation}
for some densities~$\hat{\sigma}$, and $\hat{\mu}$ (that will, of course, depend on~$\by$
from the previous expressions). These are determined by imposing the conditions
in~\eqref{continuityx2} using the plane wave representation of the free-space
Green's function:
\begin{equation}
g_k(\bx,\by) = 
\frac{1}{4\pi}\int_{-\infty}^{\infty} 
\frac{ e^{-\sqrt{\xi^2 - k^2}|x_2 - y_2|}}
{\sqrt{\xi^2 - k^2}}  e^{i(x_{1}-y_1)\xi} \, d\xi ,
\end{equation}
see~\cite{chew,sommerfeld} for details and a derivation.  After evaluating the
above expression along the interface~$x_2 = 0$, enforcing the interface
conditions on~$g^{LM}$, and some algebra, one obtains that
\begin{equation}
  \label{eq:ucdef}
  s_{+}(\bx, \by)
  = \frac{1}{4\pi}\int_{-\infty}^{\infty}
  \frac{k_{+}^2 - k_{-}^2}{\sqrt{\xi^2 - k_{+}^2}
    \left( \sqrt{\xi^2 - k_{+}^2} + \sqrt{\xi^2 - k_{+}^2}\right)^2}
   e^{-\sqrt{\xi^2 - k_{+}^2}(x_{2} + y_{2})} \, e^{i (x_{1}-y_{1}) \xi} \, d\xi \, ,
\end{equation}
where the arguments of the square roots are taken to be such
that~$\textrm{Arg}(\sqrt{z}) \in [0, \pi)$.  (A similar formula holds for
computing~$s_{-}$.)

With the Green's function in hand, we may define the layered medium single and
double layer potentials for an inclusion $\Gamma$ by
\begin{equation}
  \begin{aligned}
    \cS_{\Gamma}^{LM}[\sigma](\bx) &= \int_{\Gamma} g^{LM}(\bx,\by) \,
                                     \sigma(\by) \, ds, \\
    \cD_{\Gamma}^{LM}[\mu](\bx) &= \int_{\Gamma}
                                     \left( \bn(\by) \cdot \nabla_{\by}
                                     g^{LM}(\bx,\by)  \right) \,
                                     \mu(\by) \, ds ,
  \end{aligned}  
\end{equation}
where we assume that~$\bx$ is in the upper half space, as per the earlier
calculation.  The jump relations on~$\Gamma$ are identical to their free space
counterparts given in~\eqref{eq:jumprel}.

For $\bx$ on a second curve $\Gamma'$ with outward normal
$\bn(\bx)$, we define their normal derivatives by
\begin{equation}
  \begin{aligned}
\cS{'}_{\Gamma}^{LM}[\sigma](\bx) &= \bn(\bx) \cdot \nabla
\cS_{\Gamma}^{LM}[\sigma](\bx), \\
\cD{'}_{\Gamma}^{LM}[\mu](\bx) &= \bn(\bx) \cdot \nabla
\cD_{\Gamma}^{LM}[\mu](\bx)  .
\end{aligned}
\end{equation}

\subsection{The layered medium scattering matrix}

So long as an obstacle $\Gamma$ and its enclosure $P$ lie entirely in the upper 
half space, the scattered field
can be expressed using the combined field representation 
\begin{equation}
  \usc = \cD_{\Gamma}^{LM}[\sigma](\bx) + ik\cS_{\Gamma}^{LM}[\sigma](\bx)
  \, , \quad \bx \in \Omega^{c} \,.
\label{lmscat}
\end{equation}
Note that by using the layered medium Green's function in the above
representation,~$\usc$ automatically satisfies the interface conditions, the
proper radiation conditions for the layered media, as well as the same Green's
identities as the free space
kernel~\cite{epstein2023solving1,epstein2023solving2,epstein2024solving3}.  We
will let $\mtx{A}_{\Gamma}^{P}$ denote the scattering operator for a single
scatterer in the presence of the layered medium as well, since the context will
be clear. Following the analysis of Section~\ref{sec:scat}, we have
\begin{equation}
\begin{bmatrix}
    \usc|_{P} \\
    \dudnsc|_{P}
\end{bmatrix}
= \mtx{A}^{\Gamma}_{P} 
\begin{bmatrix}
\uin|_{P}  \\
\dudnin|_{P} 
\end{bmatrix} \, ,
\end{equation}
where
\begin{equation}
\label{eq:scatmatlm}
\mtx{A}^{\Gamma}_{P}
= 
\begin{bmatrix}
{\cD}^{LM}_{\gammatop} + ik {\cS}^{LM}_{\gammatop} \vspace*{1ex}\\
{\cD}'^{LM}_{\gammatop} + ik {\cS}'^{LM}_{\gammatop}
\end{bmatrix}
\left( \frac{1}{2} + \mtx{D}^{LM}_{\Gamma} + ik\mtx{S}^{LM}_{\Gamma} \right)^{-1}
\begin{bmatrix}
{\cD}^{LM}_{\ptogamma} & -{\cS}^{LM}_{\ptogamma}
\end{bmatrix} .
\end{equation}
In the above, we have also assumed that~$\uin$ satisfies the layered media interface conditions.

\begin{remark}
  In practice, it is often useful to take~$\uin$ as a free-space planewave, in
  which case an extra correction term, e.g. a reflection obtained from Snell's
  law, must be added to ensure that the total field satisfies those layered media
  interface conditions.
\end{remark}

The full multi-particle system still takes the form
\eqref{eq:multi-part-scat-main}, with the layered medium scattering matrix
$\mtx{A}^{\Gamma}_{P}$ used in the definition of 
$\cA$ in~\eqref{Adef}, and with
\begin{equation}
  \label{eq:tij}
\mtx{T}_{ij} = \begin{bmatrix}
{\cD}^{LM}_{P_{j} \rightarrow P_{i}} & 
-{\cS}^{LM}_{P_{j} \rightarrow P_{i}} \vspace*{1ex}\\
{\cD}'^{LM}_{P_{j} \rightarrow P_{i}} & 
-{\cS}'^{LM}_{P_{j} \rightarrow P_{i}}
\end{bmatrix} 
\end{equation}
used in the definition of $\cT$ in~\eqref{Tdef}.

\section{Discretization}
\label{sec:disc}
Up to this point, we have described the formalism for the construction of
scattering operators in terms of continuous operators acting on functions. In
practice, we need to numerically discretize the inclusion boundaries~$\Gamma_j$,
the enclosure boundaries~$P_j$, and in the layered media setup, the associated
Sommerfeld integrals which provide the corrections for the interface conditions.

\subsection{Scattering in free space}

The discretization of the obstacle itself is needed in order to determine the
density~$\sigma$ in~\eqref{eq:sig_sol}, or its analog in the layered medium
case.  There are several relatively standard approaches to this in the
integral equation literature. For this work, we choose to use the software
package \texttt{chunkie}~\cite{chunkie}. This particular numerical solver uses a
high order discretization of~$\Gamma_j$ based on adaptively determined panels,
and applies a generalized Gaussian quadrature~\cite{ggq10} to accurately
discretize the integral equations. A Nystr\"{o}m-type discretization is used.
Within the package, and depending on the problem size, the integral equation can
be solved directly, iteratively with fast multipole acceleration~\cite{fmm2d},
or with a fast direct solver \cite{FLAM}. Since it plays a ``black box'' role in
the present paper, we omit further details.

The proxy surface boundary $P$, on the other hand, and its efficient handling is
an important part of the overall procedure.  We choose to form~$P$ as a closely
fitting rectangle around~$\Gamma$, and discretize it using composite 16th-order
Gauss-Legendre panels.  Letting $n_P$ denote the total number of points
along~$P$, it is well-known that for non-oscillatory kernels -- even singular or
weakly singular ones -- it is sufficient for~$n_P$ to be of the
order~$\cO(\log{(1/\dmin)})$, where~$\dmin = \min_{i\neq j} d(P_{i}, P_{j})$ and
with~$d(P_{i}, P_{j})$ defined as the smallest distance between the proxy
surfaces. In the oscillatory regime, for scatterers which are many wavelengths
in diameter the number of points~$n_P$ must grow linearly with the size of~$P$
in order for both the outgoing and incoming fields to be sufficiently sampled.
See~\cite{cheng2006,greengard1998accelerating,bremer-quad-2015} for details.

For the multiple scattering setup with~$M$ inclusions, we let $N = M \cdot n_P$
and let~$\by_j$, for~$j = 1,\ldots,N$, denote the set of all proxy surface
discretization points, which we will refer to as {\em proxy points}.
We furthermore note that each fully discretized single particle scattering
matrix is a~$2n_P \times 2n_P$ matrix, and that the total scattered
field~\eqref{mpsfree} in the free space case is represented as the sum of all individual scattered fields from all inclusions:
\begin{equation}
\label{uscfreedisc}
\usc(\bx) = \sum_{\substack{j=1\\ \bx \neq \by_{j}}}^{N} w_j \Big( 
 \left( \bn_{j} \cdot \nabla_{\by}
  g_{k}(\bx,\by_{j}) \right) \mu_{j} - g_{k}(\bx, \by_{j}) \, \rho_j   \Big) \, ,
\end{equation}
where~$\rho_j = \frac{\partial \usc}{\partial n}(\by_j)$
and~$\mu_j = \usc(\by_j)$, and the quadrature weights along all the proxy
surfaces are given by~$w_j$.  As described early in the manuscript, the
``unknowns'' used to represent the fields are the values and normal derivatives
of the scattered field along~$P$. It is the above representation that yields
specific entries in the matrix~$\mtx{T}_{ij}$ appearing in
equation~\eqref{Tijdef}, which maps fields from each proxy surface to all the
others.

\subsection{Constructing the scattering matrix}

In order to construct a scattering matrix~$\mtx{A}_P^\Gamma$ and the
multi-particle scattering matrix~$\cA$, referring to~\cref{eq:scatmat,eq:scatmatlm}, the incoming fields are given by potentials due
to point charges and dipoles placed at the discretization points on the proxy surface $P$. The problem of scattering
from~$\Omega$ is then solved using the software library~\texttt{chunkie}, and
the outgoing scattered field is then sampled on the proxy surface~$P$ at the
same nodes~$\by_j$. For a single point charge or dipole placed at $\by_{j}$, these steps effectively provide
one column of the scattering matrix~$\mtx{A}^\Gamma_P$. 
Note that this approach is different from related procedures described in~\cite{FMPS,FMPS2D,bremer-quad-2015}, where the incident fields are typically constructed using a basis of plane waves. Using point charges, and dipoles provides additional flexibility when using commercial software to construct the scattering matrices.

\subsection{Scattering in layered media}

In addition to all of the discretization considerations mentioned in the
previous section, for layered medium problems we must add to the representation
of~$\usc$ in~\eqref{uscfreedisc} the Sommerfeld integral contribution $\usc_S$:
\begin{equation}
\usc(\bx) = \sum_{\substack{j=1\\ \bx \neq \by_{j}}}^{N} w_j \Big( 
 \left( \bn_{j} \cdot \nabla_{\by}
  g_{k}(\bx,\by_{j}) \right) \mu_{j} - g_{k}(\bx, \by_{j}) \, \rho_j   \Big) + \usc_S(\bx),
\end{equation}
where in the upper half space, we have
\begin{equation}
  \usc_S(\bx) = \sum_{\substack{j=1\\ \bx \neq \by_{j}}}^{N} w_j \Big(
  \left( \bn_{j} \cdot \nabla_{\by} s_{+}(\bx,\by_{j}) \right) \mu_{j}
  -s_{+}(\bx, \by_{j}) \rho_{j}
  \Big).
\end{equation} 
Above, as in~\eqref{uscfreedisc}, we have set
$\rho_j = \frac{\partial \usc}{\partial n}(\by_j)$ and~$\mu_j = \usc(\by_j)$.
The above representation, along with the associated one obtained from using~$s_-$
in a Sommerfeld representation in the lower half space with material
parameter~$k_-$, ensures that the fields automatically satisfy the interface
conditions along~$x_2 = 0$. 

After some algebraic rearrangement, it is straightforward to see that
\begin{equation}
\label{eq:usommcor}
\usc_S(\bx)  
= \frac{1}{4\pi}\int_{-\infty}^{\infty} 
\frac{ (k_+^2 - k_{-}^2)  e^{-\sqrt{\xi^2 - k_{+}^2} \, x_2} } 
{\sqrt{\xi^2 - k_{+}^2} 
\left( \sqrt{\xi^2 - k_{+}^2} + \sqrt{\xi^2 - k_{+}^2}\right)^2}
W(\xi) \, e^{i x_1 \xi} \, d\xi \, ,
\end{equation}
where
\begin{equation}
  W(\xi) = \sum_{j=1}^N
  \left[  \mu_j
    \left(  - i \xi n_1^j - \sqrt{\xi^2 - k_{+}^2} n_2^j \right)
    - \rho_j \right]
 e^{-\sqrt{\xi^2 - k_{+}^2}x_2^j} \, e^{-i  x_1^j \xi} ,
\end{equation}
and with~$\by_j = (y_1^j,y_2^j)$ and $\bn_j = (n_1^j,n_2^j)$.

Let us now assume that the scatterers and targets live in the upper half-space
in a box of dimension $[-A,A] \times [\delta,B]$ with $B> \delta > 0$.  Once
$\xi > k_+$, note that the term $W(\xi)$ above is exponentially decaying at a
rate of the order $e^{-\xi \delta}$ so that exponential accuracy with precision
$\varepsilon$ is achieved by truncating the integral in \eqref{eq:usommcor} at
$\xi_{max} = k_+ + \cO( \log(1/\varepsilon)/\delta)$.  Since the integrand
includes the oscillatory term $e^{i x_1 \xi}$, approximately~$\cO(A \xi_{max})$
points are needed for Nyquist-rate sampling, even if the integrand were
smooth. There are, however, singularities in~$\xi$ at both~$\pm k_+$
and~$\pm k_{-}$.

Without going into detail, there are many quadrature methods available to
address these singularities, including generalized Gaussian
quadratures~\cite{ggq10}, end-point corrected trapezoidal rules~\cite{alpert},
and adaptive Gaussian quadrature, where the panels are dyadically refined in the
vicinity of $\xi = \pm k_{\pm}$. In the present work, we use a combination of end-point corrected
trapezoidal rules and adaptive refinement.  We denote the full set of quadrature
nodes and weights by~$\xi_l$ and~$r_l$, respectively, for $l =
1,\ldots,N_\xi$. This leaves us with the need to compute
\[
W(\xi_l) = \sum_{j=1}^N  \left[ 
  \mu_j ( - i \xi_l n_1^j - \sqrt{\xi_l^2 - k_{+}^2} n_2^j )
  - \rho_j \right]
 e^{-\sqrt{\xi_l^2 - k_{+}^2}x_2^j} \, e^{-i  x_1^j \xi_l} ,
\]
and once the $W(\xi_l)$ are known, the correction becomes
\begin{equation}
\usc_S(\bx)  
= \sum_{l=1}^{N_\xi} r_l
\frac{ (k_+^2 - k_{-}^2) \,  e^{-\sqrt{\xi_l^2 - k_{+}^2} \, x_2} } 
{\sqrt{\xi_l^2 - k_{+}^2} 
\left( \sqrt{\xi_l^2 - k_{+}^2} + \sqrt{\xi_l^2 - k_{+}^2}\right)^2}
W(\xi_l) \, e^{i x_1 \xi_l } .
\end{equation}

There is a large literature devoted to accurately and efficiently computing with
Sommerfeld integrals and representations, see for
example~\cite{okhmatovski2024,lai2018new, lindell1984exact,hochman2009numerical,
  taraldsen2005complex}, and a detailed comparison of various approaches for the
computation of layered medium Green's functions will be carried out at a later
date. It is more instructive to see what the actual values of $N_\xi$ are in
specific examples. We will note, however, that as $\delta$ gets smaller, the
Sommerfeld integral becomes more and more troublesome to evaluate. This is a
known difficulty and there are many remedies, but all lead to additional
algorithmic complexities, so we assume for the sake of simplicity that
$2 \pi \delta \cdot \max(k_+,k_{-}) > 0.1$ -- that is to say, the sources and
targets are at least~$0.1$ wavelengths away from the interface.

\section{Fast algorithms for multi-particle scattering}

We assume here that the individual scattering matrices are modest in size so
that applying $\cA$ in~\eqref{Adef} can be done directly at a cost of the order
$O(M n_{p}^2)$.  For large-scale problems, it is the cost of applying $\cT$
that dominates.  However, as is clear from \eqref{uscfreedisc}, this can be done
with $O(N \log N)$ work using the fast multipole
method~\cite{zhang2007fast,FMPS2D,FMPS,lai2022fast}. The two dimensional fast
multipole method library, see~\cite{fmm2d}, can be used for such calculations.

The calculation of the values~$W(\xi_l)$ above, and the subsequent Sommerfeld
integrals, is slightly non-standard. However, it can also be carried out
in~$O(N_\xi \log N_\xi + N)$ work using the non-uniform FFT (NUFFT) and local
interpolation. See~\cite{finufftlib,finufft,barnett2021acha} for more
information regarding the NUFFT, and~\cite{FMPS2D} for a discussion regarding
the interpolation procedure. Alternative, there exist special-case FMM-based
algorithms for computing with the layered media Green's function,
see~\cite{cho2021adapting,
  cho2018heterogeneous,cho2012parallel,gurel1996fast,geng2001fast} for a more
detailed discussion.

\section{Numerical results}
\label{sec:results}

In this section we provide descriptions and results for several numerical
experiments to validate the ideas described in the earlier sections.
For the examples in~\Cref{subsec:accuracy}, the scatterers are a
collection of ellipses parametrized via
\begin{equation}
  E_{a,b}(t; \bc) = (a \cos{(t)} + c_{1}, b \sin{(t)} + c_{2}), \qquad
  t\in[0,2\pi),
\end{equation}
with~$\bc = (c_{1}, c_{2})$. For the examples
  in~\Cref{subsec:multi,subsec:multi-lm}, the scatterers are a collection of star-shaped
ellipses with parametrization
\begin{equation}
  W_{a,b}(t; \bc) = r(t) (a \cos{(t)} + c_{1}, b \sin{(t)} + c_{2})
\end{equation}
where $r(t) = 1 + 0.1\cos{(7t)}$.  The layer potentials for the construction of
the scattering matrices and the proxy surfaces are discretized using the
procedure described in Section~\ref{sec:disc}.  The unknowns are scaled by the
square roots of the smooth quadrature weights for better numerical conditioning~\cite{bremer2012corner}.

In the examples below, each proxy surface is assumed to be identical up to a
shift, and discretized with $n_{p}$ points. The variable~$\varepsilon_{a}$
denotes the accuracy in the computed scattered field at the discretized points
on the proxy surface as compared to either the exact solution computed by
directly discretizing the scatterers, or a self-convergence result obtained by a
finer discretization of the proxy surfaces.  In all cases, since the scatterers
are sufficiently small as measured in wavelengths, the scattering matrices are
computed using dense linear algebra. The surface of the scatterer is
sufficiently discretized so that scattering matrices are computed to a numerical
accuracy of $10^{-12}$, and~$\efmm$ denotes the tolerance used for applying~$\cT$ using fast multipole methods and NUFFTs.  In this work, we use the
software packages \texttt{fmm2d}~\cite{fmm2d} and
\texttt{finufft}~\cite{finufft} for the fast evaluation of $\cT$.  Finally,
$\timeit$ denotes the CPU time for applying the discretized multiple-scattering
equations, $\timetot$ denotes the total solve time, and $\nit$ denotes the number
of GMRES iterations required for the relative residual to drop below the
prescribed tolerance $\egmres$.

Unless stated otherwise, the incident field for testing the accuracy of the
solvers is an incoming plane wave in the $x_{1}$ direction for the free space
problems, i.e. $\uin = e^{i k x_{1}}$. For examples in the layered medium
geometry, the incident field is an incoming plane wave propagating in the
direction~$\boldsymbol{d} = (\cos{(\theta)}, \sin{(\theta)})$, with
$\theta = \pi/3$, i.e
\begin{equation}
  \uin_{\textrm{lm}}(\bx)
  = e^{\left(i k_{+}( \cos{(\theta)x_{1}} + \sin{(\theta)}x_{2}  )\right)} 
  - \frac{k_{+}\sin{(\theta)}+ \sqrt{k_{-}^2 - k_{+}^2 \cos^2{(\theta)}}}{k_{+} \sin{(\theta)}-\sqrt{k_{-}^2 - k_{+}^2 \cos^2{(\theta)}}} e^{\left(i k_{+}( \cos{(\theta)x_{1}} - \sin{(\theta)}x_{2}  )\right)} \, .
\end{equation}
Note that $\uin_{\textrm{lm}}$ includes the contribution from Snell's law so
that it satisfies the layered medium equation in $\mathbb{R}^{2}$.

\subsection{Effects of distance, aspect ratio, and wave number
  on number of proxy points}
\label{subsec:accuracy}

In this section, we illustrate the dependence of the number of proxy points
required to obtain a fixed tolerance, as a function of the distance between the
two scatterers, the aspect ratio of the scatterers and the wave number; where
two out of the three parameters are held constant.  In particular consider two
ellipses $\mathcal{E}_{1} = E_{a/2,1/2}(t; (0,0))$
and~$\mathcal{E}_{2} = E_{a/2,1/2}(t, (0, 1 + d))$, where the ellipses are
separated by a distance $d$ and the aspect ratio $\eta$ of the ellipses is
$a$. The proxy surfaces are rectangles with the same centers as the ellipses and
have side lengths $2a + 2d/3$ and $2 + 2d/3$, respectively. This choice of
spacing ensures that the distance between the two proxy surface rectangles is
the same as the distance between the proxy surface and the ellipse.

In~\Cref{fig:ex2llip}, we plot $\varepsilon_{a}$ as a function of $n_{p}$ for
the following three setups:
\begin{enumerate}
  \item[(a)] on the left, for $a = 10$ and $d=1$, and three
    different wave numbers $k = 2\pi, 4\pi$, and $8 \pi$;
  \item[(b)] in the middle, for $a=10$, $k= 4\pi$, and $d = 0.5, 1$ and $2$; and
  \item[(c)] on the right, for $d=1$, $k = 20\pi/a$, and $a= 10, 20$, and $40$.
\end{enumerate}
For a fixed wave number, the method is spectrally convergent and the number of
points required to achieve a fixed accuracy increases linearly in $k$. As a
function of the separation of distance between the two ellipses, the number of
points required grows like $\cO(1/d)$.  This issue can be remedied using an
adaptive discretization in the vicinity of the regions where the ellipses are
close-to-touching, reducing the required number of points to
$\cO(\log{(1/d)})$. The number of points required as a function of aspect ratio
$\eta$ also grows like $O(\eta)$, as under a rescaling of the coordinates an
increase in aspect ratio is equivalent to a reduction in the distance between
the ellipses.

We repeat the same experiment for the layered medium setup where
$k_{-} = 1.3 k_{+}$. The ellipses are given
by~$\mathcal{E}_{1} = E_{a/2,1/2}(t; (0,2))$
and~$\mathcal{E}_{2} = E_{a/2,1/2}(t, (0, 3 + d))$, i.e. the ellipse closest to
the interface is a distance of $2$ away from the
interface. In~\Cref{fig:ex2llip_lm}, we plot analogous results for the layered
medium problem, and the trends in the error as a function of the number of proxy
points are similar to the corresponding results in the free-space setup.

\afterpage{
  \clearpage
  
  \begin{figure}[!h]
    \centering
    \begin{subfigure}[t]{0.55\linewidth}
      \centering
      \includegraphics[width=1\linewidth]{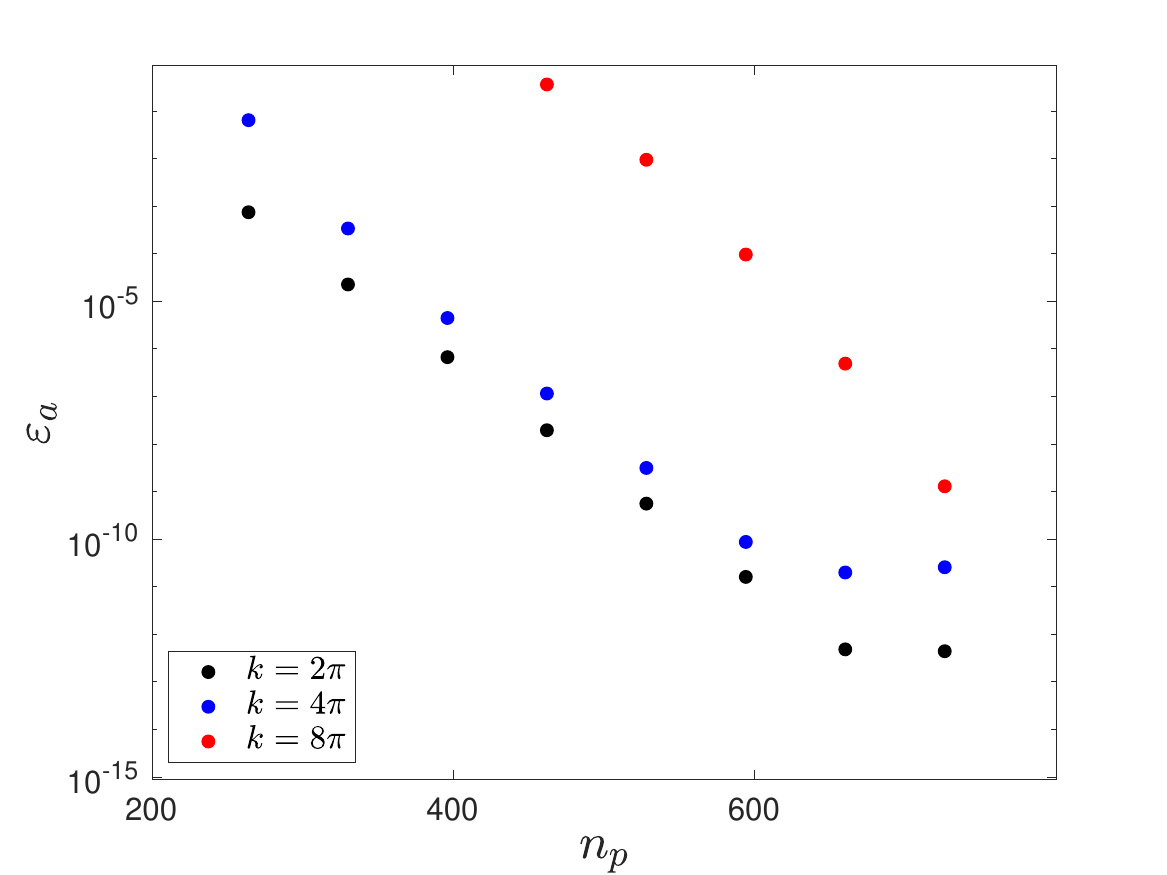}
      \caption{Dependence on $k$ for a fixed aspect
      ratio $a = 10$, and separation $d = 1$.}
    \end{subfigure} \\
    \begin{subfigure}[t]{0.55\linewidth}
      \centering
      \includegraphics[width=1\linewidth]{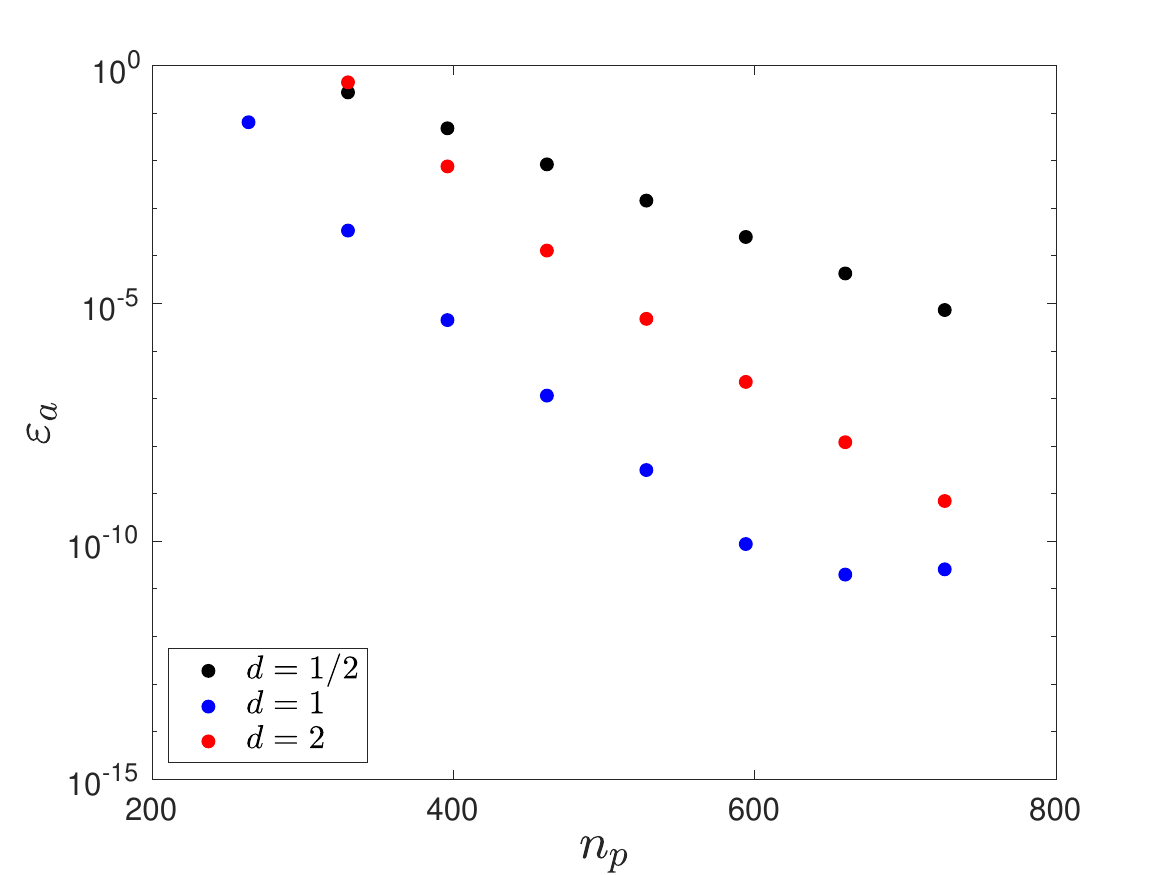}
      \caption{Dependence on $d$ for
      fixed wavenumber $k=4\pi$, and aspect ratio $a = 10$.}
    \end{subfigure} \\
    \begin{subfigure}[t]{0.55\linewidth}
      \centering
      \includegraphics[width=1\linewidth]{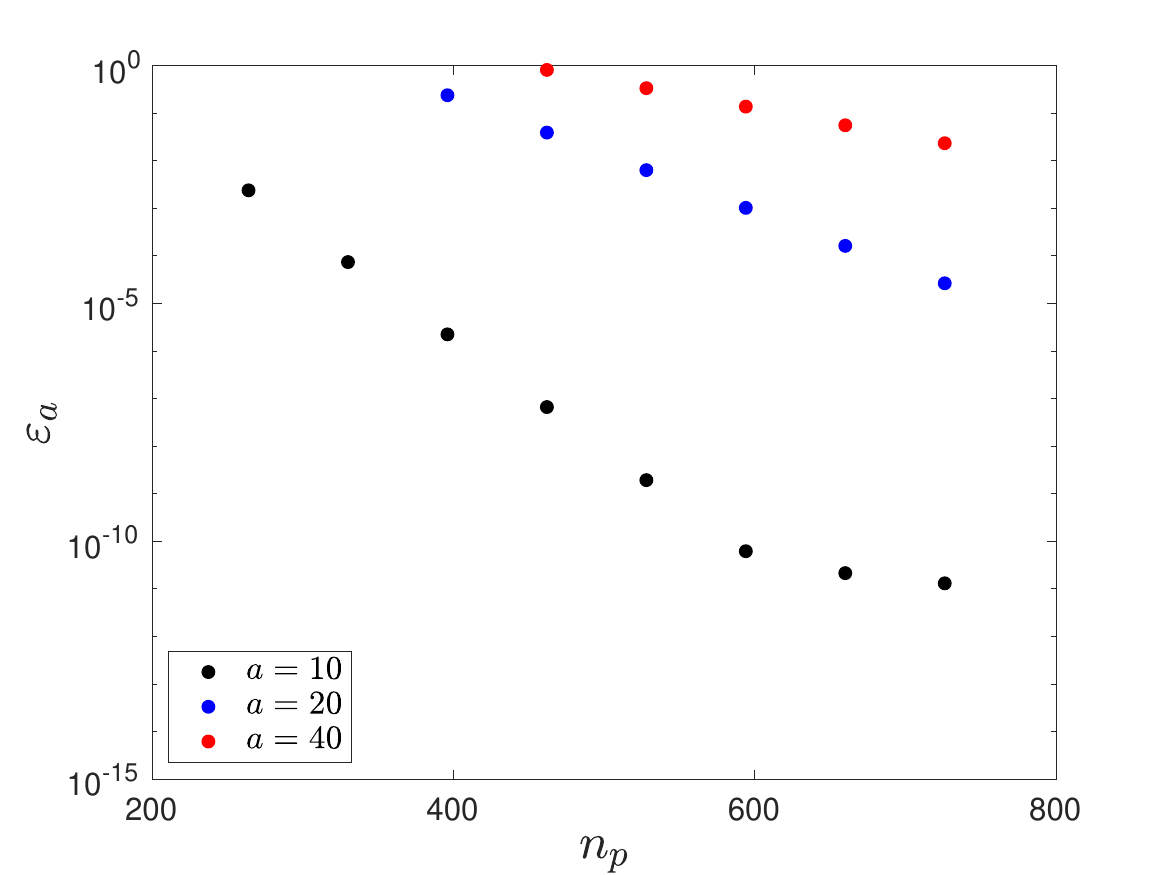}
      \caption{Dependence on the aspect ratio $a$ with wave number
        $k = 20 \pi /a$ and fixed distance $d=1$.}
    \end{subfigure}
    \caption{Shown is~$\varepsilon_{a}$ as a function of $n_{p}$ for free space
      scattering from two ellipses.}
  \label{fig:ex2llip}
\end{figure}

\clearpage
}

\afterpage{
  \clearpage
  
\begin{figure}[!h]
  \centering
  \begin{subfigure}[t]{0.55\linewidth}
    \centering
    \includegraphics[width=1\linewidth]{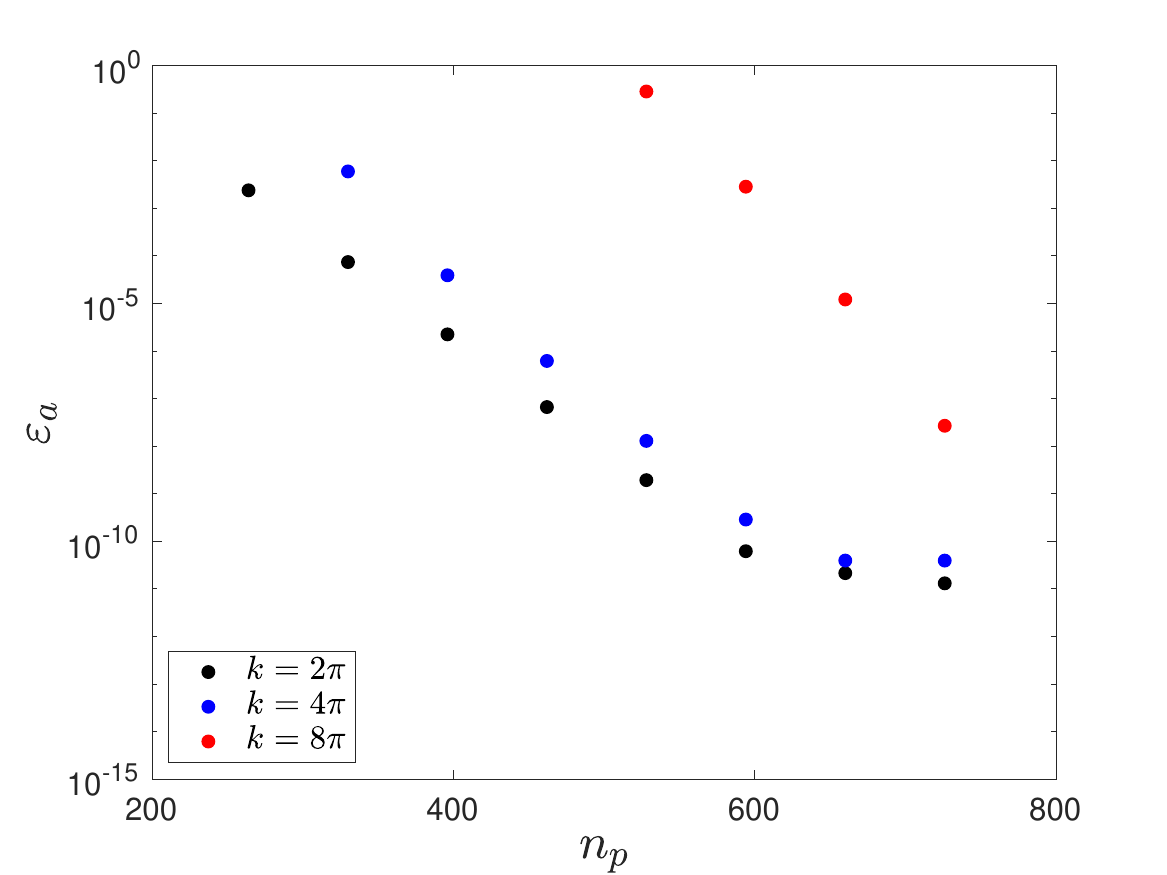}
    \caption{Dependence on $k$ for a fixed
    aspect ratio $a = 10$, and separation $d = 1$.}
  \end{subfigure} \\
  \begin{subfigure}[t]{0.55\linewidth}
    \centering
    \includegraphics[width=1\linewidth]{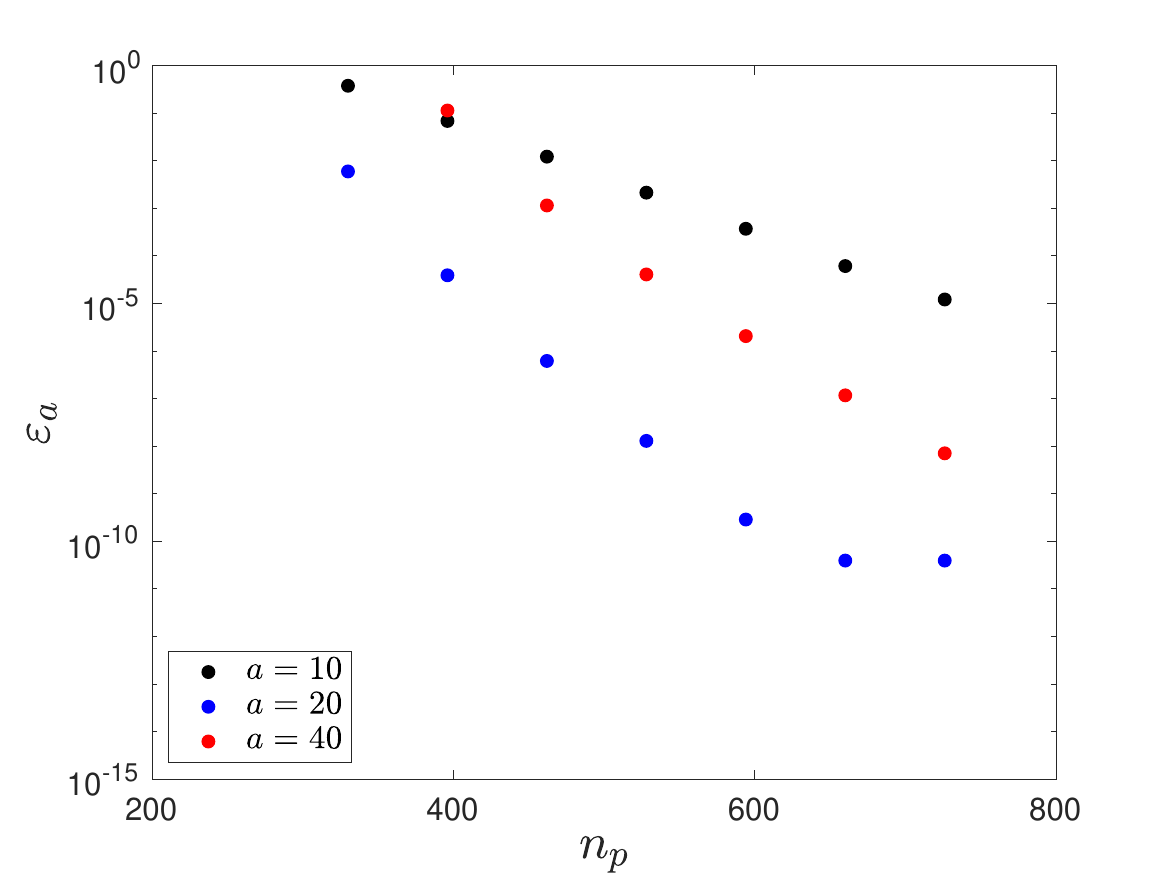}
    \caption{Dependence on $d$
    for fixed wavenumber $k=4\pi$, and aspect ratio $a = 10$.}
  \end{subfigure} \\
  \begin{subfigure}[t]{0.55\linewidth}
    \centering
    \includegraphics[width=1\linewidth]{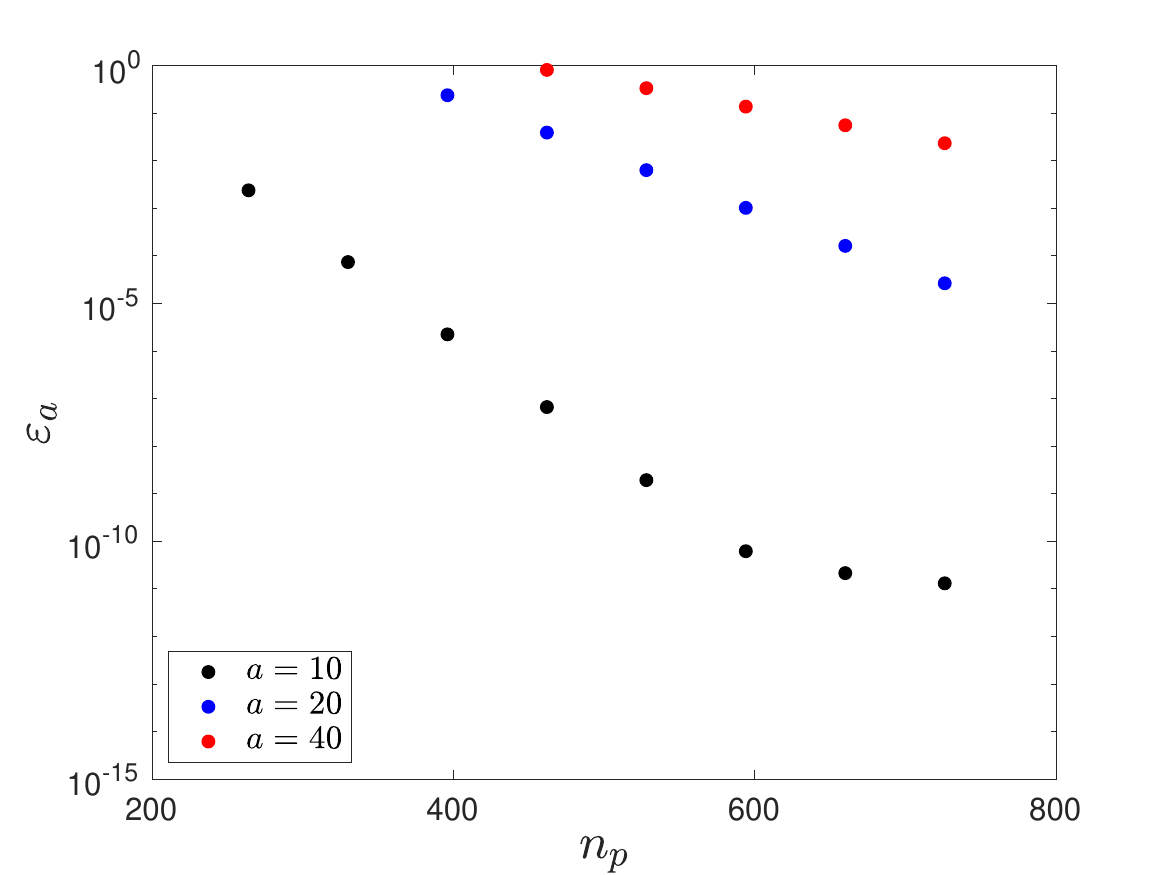}
    \caption{Dependence on the aspect ratio $a$ with wave number $k = 20 \pi /a$ and
    fixed distance $d=1$.}
  \end{subfigure}
  \caption{Shown $\varepsilon_{a}$ as a function of $n_{p}$ for scattering from
    two ellipses in a layered medium.}
  \label{fig:ex2llip_lm}
\end{figure}

\clearpage
}

\subsection{Multi-particle scattering}
\label{subsec:multi}

In this example, we consider multiple particle scattering from a $41\times 21$
lattice with a waveguide channel removed, see the dark objects in
Figure~\ref{fig:ex_wellip_lattice}. This configuration is related to the the
photonic crystal example in~\cite{gillman2015spectrally}, and we refer to this
lattice as the photonic crystal. Each scatterer is a star-shaped ellipse
$W_{a,b}(t; \bc_{ij})$, with $a=0.05/3$ and $b=0.1/3$, and the centers of the
ellipses are given by
\begin{equation*}
  \bc_{ij} = \begin{cases}
    \left(-1+(i-1)0.05,-1+(j-1)0.1\right), & \text{$j$ even},\\
    \left(-1+(i-1)0.05,-0.95+(j-1)0.1\right), & \text{otherwise},
  \end{cases}    
\end{equation*}
for $i=1,\ldots,41$ and $j=1,2\ldots,21$.  After removing the waveguide channel,
there are a total of $827$ scatterers in the photonic crystal.  We set
$k= 60 \pi$, which corresponds to the crystal being approximately $60$
wavelengths in diameter.  To construct the scattering matrix, each scatterer was
discretized with $1792$ points, and each proxy rectangle is discretized with
$n_{p} = 180$ points with $60$ points along the vertical edges, and $30$ points
along the horizontal edges. The proxy rectangles are placed such that the
distance between the closest rectangles is equal to the distance between the
proxy surface and its corresponding star-shaped ellipse, see the red rectangles
in Figure~\ref{fig:ex_wellip_rect_lattice}.

The results were computed on an 8-core Intel Core i9 laptop with 64 GB of
memory. The CPU time for applying the discretized multiple scattering operator
was $\timeit = 0.71$s, with $\efmm = 10^{-9}$. The GMRES iterations converged to
a relative residual of $\egmres=10^{-9}$ after $\nit = 3435$ iterations, with
$\timetot = 1.1 \times 10^{4}$s. Note that the large discrepancy between the
$\timeit \cdot \nit$ and $\timetot$ can be accounted for by the time taken in
memory movement between RAM and L1-cache. The solution is computed to a relative
accuracy of $1.1 \times 10^{-6}$ estimated via a self-convergence test with a
reference solution computed using $n_{p} = 360$. Furthermore, note that the
solution in the interior of the proxy surfaces is identically zero due to a
version of Green's identity. The field in the interior of the proxy surfaces and
exterior of scatterers can be computed using a different Green's identity.  We
present the real part of the scattered and the absolute value of the total field
in Figure~\ref{fig:ex_m_ellip}.

\afterpage{
  \clearpage
  
  \begin{figure}[!h]
    \centering
    \begin{subfigure}[b]{0.75\linewidth}
      \centering
      \includegraphics[width=1\linewidth]{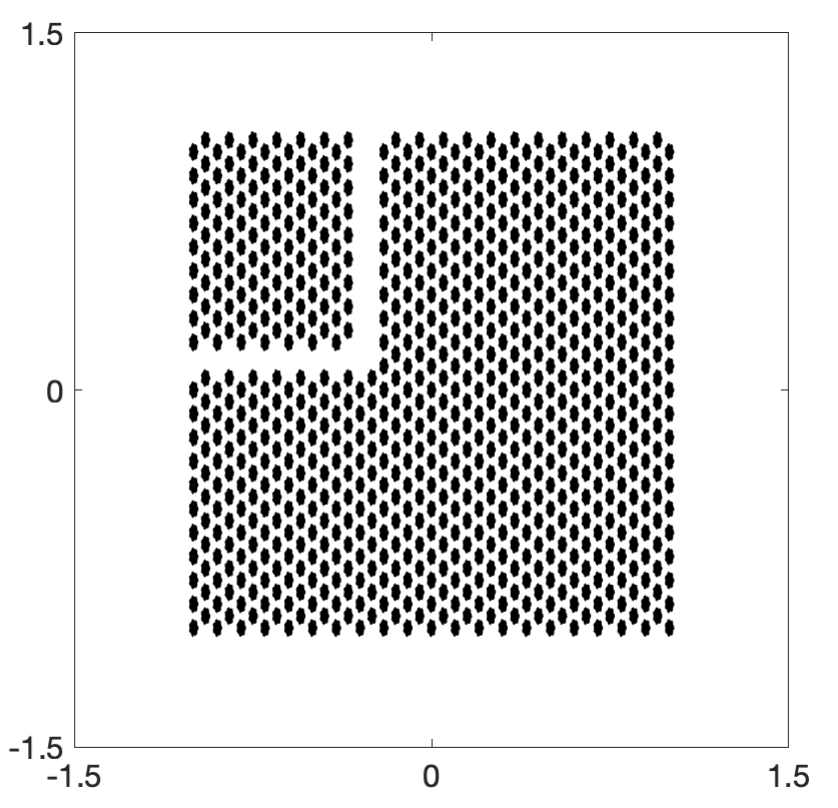}
      \caption{Photonic crystal.}\label{fig:ex_wellip_lattice}
    \end{subfigure}\\
    \begin{subfigure}[b]{0.75\linewidth}
    \centering
    \includegraphics[width=1\linewidth]{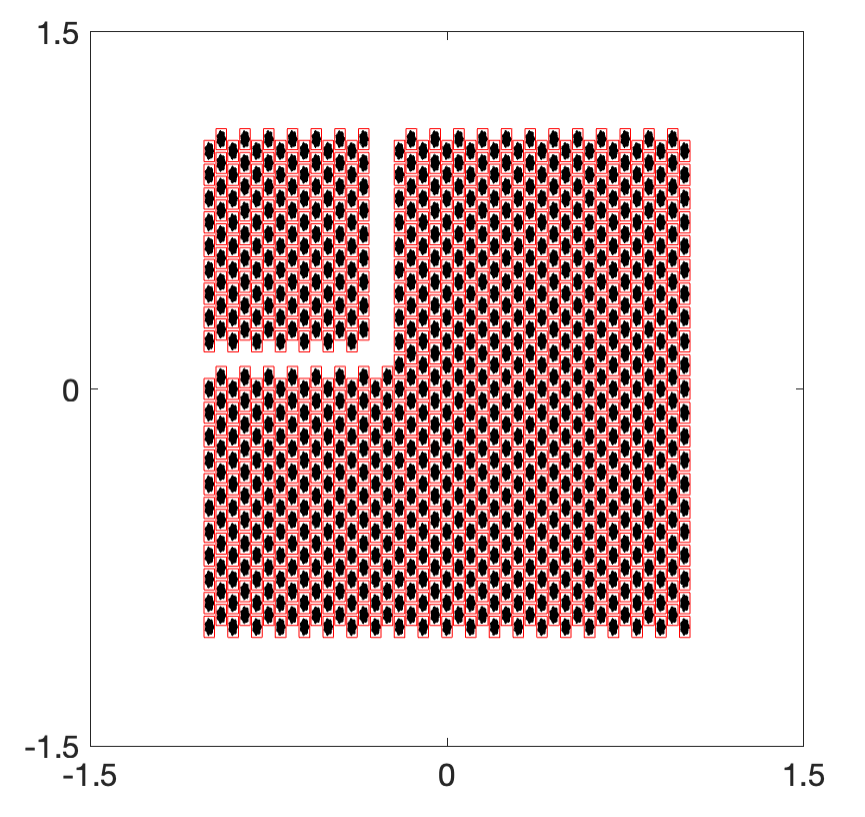}
    \caption{Photonic crystal with the proxy surfaces.}
    \label{fig:ex_wellip_rect_lattice}
  \end{subfigure}
  \caption{The geometry and proxy surfaces in the photonic crystal setup.}
  \label{fig:photonic_geometry}
\end{figure}

\clearpage
}

\afterpage{
  \clearpage
  
  \begin{figure}[!h]
    \centering
    \begin{subfigure}[b]{0.75\linewidth}
      \centering
      \includegraphics[width=1\linewidth]{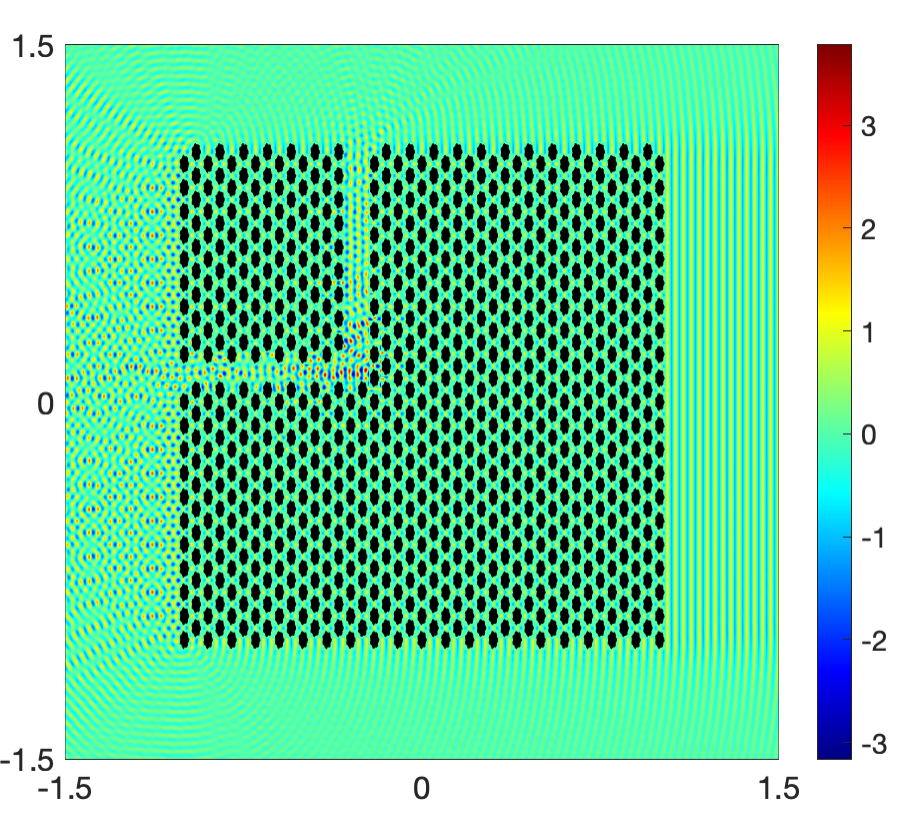}
      \caption{Real part of $\usc$.}
      \label{fig:field_photonic_real}
    \end{subfigure} \\
    \begin{subfigure}[b]{0.75\linewidth}
      \centering
      \includegraphics[width=1\linewidth]{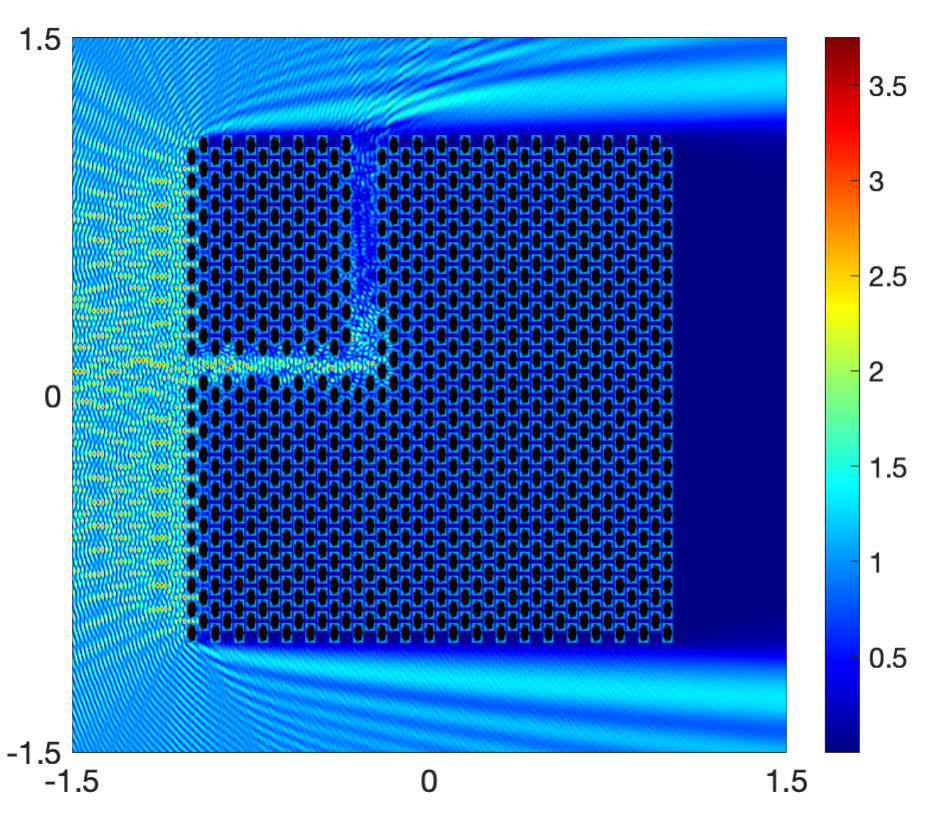}
      \caption{Absolute value of total field $|u|$.}
      \label{fig:field_photonic_abs}
    \end{subfigure}
    \caption{Multi-particle scattering from a photonic crystal.}
    \label{fig:ex_m_ellip}
\end{figure}

\clearpage
}

\subsection{Multi-particle scattering in layered media}
\label{subsec:multi-lm}
We now consider multiple scattering from a $21\times 2$ perturbed lattice of scatterers in a layered medium, see dark objects in Figure~\ref{fig:ex_wellip_lattice_lm}. This example represents a simplified model for studying cross-talk in an antenna array for the Hydrogen Intensity and Real-time Analysis experiment (HIRAX)~\cite{newburgh2016hirax, crichton2022hydrogen, kuhn2022antenna}. Each scatterer is a star-shaped ellipse $W_{a,b}(t; \bc_{ij})$ with $a = 1$, and $b=0.5$ with the centers of the ellipses $\bc_{ij}$ are given by
\begin{equation*}
  \bc_{ij} = \begin{cases}
    \left((-10 + (i-1))3 + \eta_{1,ij}, 1.6 + \eta_{2,ij} \right), & \text{j = 1, i=1,2,\ldots 21},\\
    \left((-9.5 + (i-1))3 + \eta_{1,ij}, 3.6 + \eta_{2,ij} \right), & \text{j=2, i=1,2,\ldots 20},
  \end{cases}    
\end{equation*}
with $\eta_{1,ij}$, and $\eta_{2,ij}$ are uniform random numbers in $[-0.1, 0.1]$. There are a total of $41$ scatterers in this configuration. We set $(k_{+}, k_{-}) = (\pi, 1.3\pi)$, which corresponds to the domain being approximately 40 wavelengths in length and 2.5 wavelengths in height. The number of points on the scatterer, the relative location of the proxy surfaces, and the number of points on the proxy surface are the same as in~\Cref{subsec:multi}.
\afterpage{
  \begin{figure}[!h]
      \centering
      \includegraphics[width=1\linewidth]{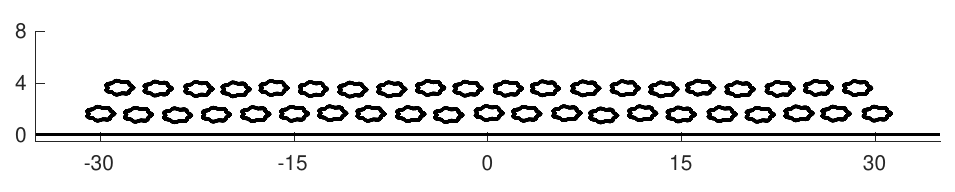}
  \caption{The geometry of scatterers in a layered medium.}
  \label{fig:ex_wellip_lattice_lm}
\end{figure}
}

The CPU time for applying the discretized multiple scattering operator
was $\timeit = 7.9$s, with $\efmm = 10^{-6}$, with the computation of Sommerfeld integral part of the Green's function being the dominant cost. The GMRES iterations converged to
a relative residual of $\egmres=10^{-6}$ after $\nit = 107$ iterations, with
$\timetot = 9 \times 10^{2}$s. The solution is computed to a relative
accuracy of $2 \times 10^{-8}$ estimated via a self-convergence test with a
reference solution computed using $n_{p} = 360$.  We
present the real part and the absolute value of the total field in the exterior of the
proxy rectangles in Figure~\ref{fig:ex_m_ellip_lm}.

\afterpage{
  
  \begin{figure}[!h]
    \centering
    \begin{subfigure}[b]{0.9\linewidth}
      \centering
      \includegraphics[width=1\linewidth]{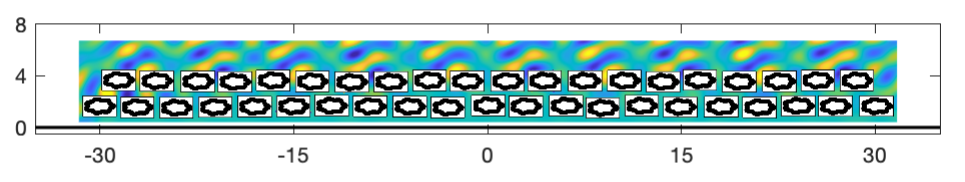}
      \caption{Real part of $u$.}
      \label{fig:field_photonic_real_lm}
    \end{subfigure} \\
    \begin{subfigure}[b]{0.9\linewidth}
      \centering
      \includegraphics[width=1\linewidth]{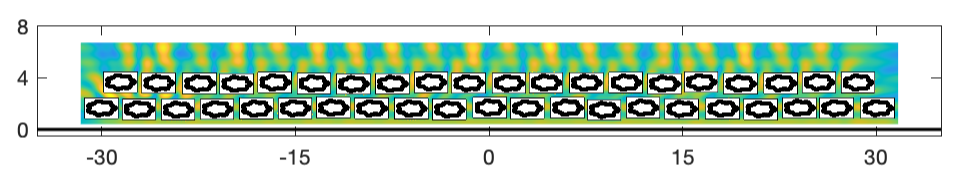}
      \caption{Absolute value of total field $|u|$.}
      \label{fig:field_photonic_abs_lm}
    \end{subfigure}
    \caption{Multi-particle scattering an array of scatterers in a layered medium.}
    \label{fig:ex_m_ellip_lm}
\end{figure}

\clearpage
}

\section{Conclusions}
\label{sec:conclusions}

In this paper we have described a general-purpose method for numerically
constructing scattering matrices for arbitrary obstacles, in particular those
that are elongated and for which standard methods would suffer from loss of
accuracy (or increased computational cost). The method is based purely on
Green's representation identities, and is compatible with black-box PDE solvers
in the sense that only values and gradients of incoming fields and outgoing
solutions are needed. The approach can also be directly applied to scattering in
layered media. Only the two-layer case was discussed in this work, but the
multi-layer case is analogous, except that each layer has a different Sommerfeld representation and additional continuity conditions, as described in \cite{FMPS2D}.

The approach used in the scattering matrix construction of this paper differs
slightly from the one used by Bremer, Gillman, and Martinsson
in~\cite{bremer-quad-2015}. In particular, the scattering matrices constructed
in that work are ones that map a fictitious density on the proxy surface to
outgoing fields on the same proxy surface. The fictitious density is responsible
for representing the incoming fields, instead of Green's representation formula
in our case. Using such an approach, as in~\cite{bremer-quad-2015}, is
particularly useful in the construction of fast direct solvers due to the way in
which blocks of the matrix are compressed -- these blocks naturally map sources
to potentials. It is likely that the approach using Green's identities of this
work could also be used in the construction of general purpose fast direct
solvers, and investigating this is ongoing work.

Lastly, it is again worth highlighting that due to the reduction of unknowns in
a multi-particle scattering problem when using scattering matrices, optimization
problems which would otherwise be computationally impossible become
tractable. Coupling the methods of this paper with such schemes is also ongoing
work.

\section*{Acknowledgments}
We would like to thank Alex Barnett, Charles Epstein, and Jeremy Hoskins for many useful discussions.
The work of C.~Borges was supported in part by the Office of Naval Research under award number N00014-21-1-2389. The Flatiron Institute is a division of the Simons Foundation.

\section*{Data Availability}
Not applicable.

\section*{Declarations}

\subsection*{Conflict of interest} 
The authors have no other relevant financial or non-financial interests to disclose.

\bibliographystyle{plainnat}
\bibliography{refs}
 
\end{document}